\documentclass[12pt]{article}
\usepackage{graphicx}
\usepackage{amsmath,amsthm,amssymb,enumerate}
\usepackage{euscript,mathrsfs}
\usepackage[left=2cm,right=2cm,top=3.5cm,bottom=3.5cm]{geometry}
\usepackage{color}

\usepackage{soul}

\catcode`\@=11 \@addtoreset{equation}{section}

\catcode`\@=12

\allowdisplaybreaks

\newtheorem{Theorem}{Theorem}[section]
\newtheorem{Proposition}[Theorem]{Proposition}
\newtheorem{Lemma}[Theorem]{Lemma}
\newtheorem{Corollary}[Theorem]{Corollary}

\theoremstyle{definition}
\newtheorem{Definition}[Theorem]{Definition}

\newtheorem{Remark}[Theorem]{Remark}

\newcommand{\bTheorem}[1]{
\begin{Theorem} \label{T#1} }
\newcommand{\eT}{\end{Theorem}}

\newcommand{\bProposition}[1]{
\begin{Proposition} \label{P#1}}
\newcommand{\eP}{\end{Proposition}}

\newcommand{\bLemma}[1]{
\begin{Lemma} \label{L#1} }
\newcommand{\eL}{\end{Lemma}}

\newcommand{\bCorollary}[1]{
\begin{Corollary} \label{C#1} }
\newcommand{\eC}{\end{Corollary}}

\newcommand{\bRemark}[1]{
\begin{Remark} \label{R#1} }
\newcommand{\eR}{\end{Remark}}

\newcommand{\bDefinition}[1]{
\begin{Definition} \label{D#1} }
\newcommand{\eD}{\end{Definition}}

\newcommand{\vrB}{\vr_b}

\newcommand{\Ds}{\mathbb{D}_x}

\newcommand{\vuB}{\vc{u}_b}

\newcommand{\bfphi}{\boldsymbol{\varphi}}

\newcommand{\bFormula}[1]{
\begin{equation} \label{#1}}
\newcommand{\eF}{\end{equation}}

\newcommand{\Ov}[1]{\overline{#1}}

\newcommand{\aleq}{\stackrel{<}{\sim}}

\newcommand{\ageq}{\stackrel{>}{\sim}}

\newcommand{\vr}{\varrho}

\newcommand{\vu}{\vc{u}}
\newcommand{\vm}{\vc{m}}

\newcommand{\vc}[1]{{\bf #1}}

\newcommand{\Div}{{\rm div}_x}
\newcommand{\Grad}{\nabla_x}

\newcommand{\dx}{\,{\rm d} {x}}

\newcommand{\dt}{\,{\rm d} t }

\newcommand{\intO}[1]{\int_{\Omega} #1 \ \dx}

\newcommand{\D}{{\rm d}}

\newcommand{\ep}{\varepsilon}

\def\softd{{\leavevmode\setbox1=\hbox{d}%
          \hbox to 1.05\wd1{d\kern-0.4ex{\char039}\hss}}}
\definecolor{Cgrey}{rgb}{0.85,0.85,0.85}
\definecolor{Cblue}{rgb}{0.50,0.85,0.85}
\definecolor{Cred}{rgb}{1,0,0}
\definecolor{fancy}{rgb}{0.10,0.85,0.10}

\newcommand\Cbox[2]{%
    \newbox\contentbox%
    \newbox\bkgdbox%
    \setbox\contentbox\hbox to \hsize{%
        \vtop{
            \kern\columnsep
            \hbox to \hsize{%
                \kern\columnsep%
                \advance\hsize by -2\columnsep%
                \setlength{\textwidth}{\hsize}%
                \vbox{
                    \parskip=\baselineskip
                    \parindent=0bp
                    #2
                }%
                \kern\columnsep%
            }%
            \kern\columnsep%
        }%
    }%
    \setbox\bkgdbox\vbox{
        \color{#1}
        \hrule width  \wd\contentbox %
               height \ht\contentbox %
               depth  \dp\contentbox
        \color{black}
    }%
    \wd\bkgdbox=0bp%
    \vbox{\hbox to \hsize{\box\bkgdbox\box\contentbox}}%
    \vskip\baselineskip%
}

\date{}

\begin{document}


\title{Globally bounded trajectories for the barotropic Navier--Stokes system with general boundary conditions}

\author{Jan B\v rezina \and Eduard Feireisl
\thanks{The work of E.F. was partially supported by the
Czech Sciences Foundation (GA\v CR), Grant Agreement
18--05974S. The Institute of Mathematics of the Academy of Sciences of
the Czech Republic is supported by RVO:67985840.}
\and Anton\' \i n Novotn\' y { }
}


\maketitle

\centerline{Faculty of Arts and Science, Kyushu University;}
\centerline{744 Motooka, Nishi-ku, Fukuoka, 819-0395, Japan}
\centerline{brezina@artsci.kyushu-u.ac.jp}

\centerline{and}

\centerline{Institute of Mathematics of the Academy of Sciences of the Czech Republic;}
\centerline{\v Zitn\' a 25, CZ-115 67 Praha 1, Czech Republic}

\centerline{Institute of Mathematics, Technische Universit\"{a}t Berlin,}
\centerline{Stra{\ss}e des 17. Juni 136, 10623 Berlin, Germany}
\centerline{feireisl@math.cas.cz}

\centerline{and}

\centerline{IMATH, EA 2134, Universit\'e de Toulon,}
\centerline{BP 20132, 83957 La Garde, France}
\centerline{novotny@univ-tln.fr}

\begin{abstract}

We consider the barotropic Navier--Stokes sytem describing the motion of a viscous compressible fluid interacting with the outer world through general in/out flux boundary conditions. We consider a hard--sphere type pressure EOS and show that all trajectories eventually enter a bounded absorbing set. In particular, the associated $\omega-$limit sets are compact and support a stationary statistical solution.

\end{abstract}


{\bf MSC:}
\bigskip


\section{Introduction}
\label{I}

A rigorous justification of \emph{ergodic hypothesis} in statistical physics is a fundamental and largely open problem 
in mathematical fluid dynamics closely related to understanding of turbulence. The Navier--Stokes system is a well accepted and widely 
used model to study these phenomena at the theoretical level. We consider its variant describing the motion of the 
density $\vr = \vr(t,x)$ and the velocity $\vu = \vu(t,x)$ of a compressible viscous barotropic fluid: 

\begin{equation} \label{i1}
\begin{split}
\partial_t \vr + \Div (\vr \vu) &= 0,\\
\partial_t (\vr \vu) + \Div (\vr \vu \otimes \vu) + \Grad p(\vr) &=
\Div \mathbb{S}(\Ds \vu) + \vr \vc{g},
\end{split}
\end{equation}
with the Newtonian viscous stress tensor
\begin{equation} \label{i1a}
\mathbb{S}(\Ds \vu) = \mu \left( \Grad \vu + \Grad^t \vu - \frac{2}{d} \Div \vu \mathbb{I} \right) +
\lambda \Div \vu \mathbb{I},\ \mu > 0, \ \lambda \geq 0,\ \ \Ds \vu \equiv \frac{1}{2} \Big( \Grad \vu + \Grad^t \vu \Big).
\end{equation}

In mathematics, it is customary to supplement \eqref{i1}, \eqref{i1a} with a kind of conservative boundary conditions, among which the most popular the no--slip 
\[
\vu|_{\partial \Omega} = 0,
\]
where $\Omega \subset R^d$ is the physical domain occupied by the fluid. The effect of the ``outer world'' is therefore 
represented by the driving force $\vc{g}$ that can be taken random in certain models of turbulence, cf. Yakhot and  Orszag  \cite{YakOrs}. 
A proper choice of $\vc{g}$ should incorporate in a ``statistically equivalent manner'' the influence of both initial and boundary conditions on the fluid motion, see \cite{YakOrs}. Although this approach has been applied successfully in the context of \emph{incompressible} fluids, see e.g. the monograph of Kuksin and Shyrikian \cite{KukShi} and the references therein, the more complex models 
of compressible fluids call for a refined treatment. Motivated by \cite{FanFeiHof}, we consider the physically relevant in/out flow boundary conditions  
\begin{equation} \label{i2}
\vu|_{\partial \Omega} = \vu_b,
\end{equation}
where $\vu_b$ is a given velocity field. Accordingly, we decompose the boundary
\[
\partial \Omega = \Gamma_{\rm in} \cup \Gamma_{\rm out},
\ \Gamma_{\rm in} = \left\{ x \in \partial \Omega \ \Big|\
 \vu_B(x) \cdot \vc{n}(x) < 0 \right\}
\]
and prescribe the density (pressure),
\begin{equation} \label{i3}
\vr|_{\Gamma_{\rm in}} = \vr_b.
\end{equation}
The (modulated) total energy of the fluid reads
\begin{equation} \label{i3a}
\intO{ E\left(\vr, \vu \Big| \vu_b \right) },\   E\left(\vr, \vu \Big| \vu_b \right)
\equiv \left[ \frac{1}{2} \vr |\vu - \vu_b|^2 + P(\vr)   \right],
\end{equation}
where $P$ is the pressure potential, 
\begin{equation} \label{i3b}
P'(\vr) \vr - P(\vr) = p(\vr).
\end{equation}
In most of the real world applications, the driving force $\vc{g}$ is the gradient of the gravitational potential 
$\vc{g} = \Grad G$ that can be incorporated in the total energy as 
\[
- \intO{ \vr G }.
\]

We are interested in fluid flows with globally bounded energy, 
\begin{equation} \label{G1}
\limsup_{t \to \infty} \intO{ E\left(\vr, \vu \Big| \vu_b \right) } \leq \mathcal{E}_\infty < \infty.
\end{equation}
As shown in \cite{FanFeiHof}, any bounded energy trajectory generates an $\omega-$limit set in the space of \emph{entire trajectories} 
(defined for $t \in R$) on which the Navier--Stokes system admits a stationary statistical solution. A stationary statistical solution 
is a stationary random process 
\[
t \in R \mapsto (\vr(t, \cdot), \vm (t, \cdot)),\ \vm \equiv \vr \vu,
\]
ranging in a suitable phase space and solving the problem \eqref{i1}--\eqref{i3} almost surely.

The existence of globally bounded trajectories for the Navier--Stokes system \eqref{i1}, with the no--slip boundary conditions, 
a general (non--potential) driving force $\vc{g}$, 
and the isentropic pressure law $p(\vr) = a \vr^\gamma$, was first established in \cite{FP14} under rather restrictive assumptions 
concerning the adiabatic exponent $\gamma$. These were partially relaxed in several papers by 
Guo, Jiang and Yin \cite{GuoJiaYin}, Jiang and Tan \cite{JiaTan},
Wang and Wang \cite{WanWan}.

Boundedness of fluid flows driven by general in/out flux inputs is much more delicate. It is actually closely related to the existence 
of (deterministic) stationary solutions to the Navier--Stokes system that is largely open for the isentropic pressure law 
$p(\vr) = a \vr^\gamma$. Following \cite{FeiNov2017} we consider the hard--sphere pressure EOS, 
\begin{equation} \label{HS}
p \in C^1[0, \Ov{\vr}) ,\ 0 < \Ov{\vr} < \infty,\ 
p'(\vr) > 0 \ \mbox{for}\ 0 < \vr < \Ov{\vr}, \ \lim_{\vr \to \Ov{\vr}-} p(\vr) = \infty.
\end{equation}
The relevant existence theory in the framework of finite energy weak solutions has been developed recently in \cite{ChoNovYan}. 
In particular, the existence of global--in--time solutions is established under the condition that 
\begin{equation} \label{i9}
\int_{\partial \Omega} \vu_b \cdot \vc{n} \ \D S_x > 0.
\end{equation}
Our goal is to show that under these circumstances the total energy of the global in time solution remains bounded as stated in 
\eqref{G1}. Moreover, we show that the constant $\mathcal{E}_\infty$ is universal, meaning the same for all trajectories. This can be rephrased as that the associated dynamical system admits a bounded absorbing set.

Our approach relies on the variant of de Rham--Hodge--Kodaira decomposition of general vector fields due to 
Kozono and Yanagisawa \cite{KozYan}, \cite{KozYan2} applied to a suitable extension of $\vu_b$ inside $\Omega$. This enables us to 
show that the driving terms in the associated energy balance can be controlled by dissipation. In general, the size of the aborbing set 
is finite but increases with increasing Reynolds number. Note that the problem is more delicate than for the stationary case studied in
\cite{FeiNov2017} as $\Div (\vr \vu) \ne 0$. 

The paper is organized as follows. In Section \ref{M}, we recall the concept of weak solution, formulate the principal hypotheses, 
and state our main result. Section \ref{B} collects all available energy and pressure estimates. A suitable extension of the boundary vector field $\vu_b$ is constructed in Section \ref{E}.  The existence of the bounded absorbing set is shown in Section \ref{S}. The paper is concluded by a short discussion on possible extensions of the main result in Section \ref{C}.

\section{Preliminaries, main results}
\label{M}

To avoid technical difficulties, we suppose that $\Omega \subset R^d$, $d=2,3$ is a bounded domain of class $C^3$, in particular that outer normal vector $\vc{n}$ is well defined on $\partial \Omega$. Possible relaxations of smoothness of $\partial \Omega$ will discussed 
in Section \ref{C}. Similarly, we suppose that the boundary data $\vu_b$, $\vr_b$ 
are restrictions of smooth vector fields $\vu_b \in C^2_c(R^d; R^d)$, $\vr_b \in C^1_c(R^d)$ to $\partial \Omega$, respectively. 
We also suppose that the pressure is given by EOS \eqref{HS}. Possible relaxation of these hypotheses will be discussed 
in Section \ref{C}.

\subsection{Global--in--time weak solutions}

\begin{Definition}[Global weak solution] \label{MD1}

We say that $(\vr, \vu)$ is a \emph{global weak solution} of the Navier--Stokes system \eqref{i1}--\eqref{i3}, with the pressure EOS 
\eqref{HS}, in $(0, \infty) \times \Omega$ if:

\begin{itemize}
\item (regularity)
\[
\begin{split}
0 \leq \vr &\leq \Ov{\vr} \ \mbox{a.a. in}\ (0,T) \times \Omega,\ 
\vr \in C_{{\rm weak,loc}}((0, \infty); L^q (\Omega)) \ \mbox{for any}\ 1 \leq q < \infty,\\
\vm \equiv \vr \vu &\in C_{{\rm weak,loc}}([0, \infty); L^{\frac{2 \gamma}{\gamma + 1}}(\Omega; R^d)),\\
(\vu - \vuB) &\in L^2_{\rm loc}((0,\infty); W^{1,2}_0(\Omega; R^d)), \\ 
p(\vr) &\in L^1_{\rm loc}((0, \infty) \times \Ov{\Omega});
\end{split}
\]

\item (equation of continuity)
\begin{equation} \label{M1}
\int_0^\infty \intO{ \Big[ \vr \partial_t \varphi + \vr \vu \cdot \Grad \varphi \Big] } \dt
= 
\int_0^\infty \int_{\Gamma_{\rm in}} \varphi \vrB \vuB \cdot \vc{n} \ \D \ S_x
\end{equation}
holds for any test function
$\varphi \in C^1_c((0,\infty) \times ({\Omega} \cup \Gamma_{\rm in}))$;

\item (momentum equation)
\begin{equation} \label{M2}
\int_0^\infty \intO{ \Big[ \vr \vu \cdot \partial_t \bfphi + \vr \vu \otimes \vu : \Grad \bfphi
+ p(\vr) \Div \bfphi \Big] } \dt = \int_0^\infty \intO{ \Big[ \mathbb{S}(\Ds \vu) : \Grad \bfphi  
- \vr \vc{g} \cdot \bfphi   \Big] } \dt
\end{equation}
holds for any test function
$\bfphi \in C^1_c((0,\infty) \times {\Omega}; R^d)$;

\item (energy inequality) 
\begin{equation} \label{M3}
\begin{split}
&- \int_0^\infty \partial_t \psi \intO{\left[ \frac{1}{2} \vr |\vu - \vuB|^2 + P(\vr) \right] } \dt  +
\int_0^\infty \psi \intO{ \mathbb{S}(\Ds \vu) : \Ds \vu } \dt \\
&  +
\int_0^{\infty} \psi  \int_{\Gamma_{\rm in}} P(\vr_b)  \vuB \cdot \vc{n} \ \D S_x \dt
\\	
&\leq -
\int_0^\infty \psi \intO{ \left[ \vr \vu \otimes \vu + p(\vr) \mathbb{I} \right]  :  \Grad \vuB } \dt + \frac{1}{2} \int_0^{\infty}
\psi \intO{ {\vr} \vu \cdot \Grad |\vuB|^2   }
\dt\\ &+ \int_0^{\infty} \psi \intO{ \mathbb{S}(\Ds \vu) : \Ds \vuB } \dt + \int_0^\infty \psi \intO{ \vr \vc{g} \cdot (\vu - \vu_b) }\dt
\end{split}
\end{equation}
holds for any test function $\psi \in C^1_c(0, \infty)$, $\psi \geq 0$.

\end{itemize}

\end{Definition}

The existence of global--in--time weak solutions for the Navier--Stokes system \eqref{i1}--\eqref{i3} 
in the framework introduced by Lions \cite{LI4}  
for any finite energy initial data was proved in \cite{ChoNovYan} under certain additional hypotheses imposed on the pressure EOS.

As a matter of fact, the existence proof in \cite{ChoNovYan} asserts only the ``integrated'' version of the energy inequality \eqref{M3}, namely 
\[
\begin{split}
&\left[ \intO{\left[ \frac{1}{2} \vr |\vu - \vuB|^2 + P(\vr) \right] } \right]_{t=0}^\tau  +
\int_0^\tau \intO{ \mathbb{S}(\Ds \vu) : \Ds \vu } \dt \\
&  +
\int_0^\tau \int_{\Gamma_{\rm in}} P(\vr_b)  \vuB \cdot \vc{n} \ \D S_x \dt
\\	
&\leq -
\int_0^\tau \intO{ \left[ \vr \vu \otimes \vu + p(\vr) \mathbb{I} \right]  :  \Grad \vuB } \dt + \frac{1}{2} \int_0^\tau
 \intO{ {\vr} \vu \cdot \Grad |\vuB|^2   }
\dt\\ &+ \int_0^\tau \intO{ \mathbb{S}(\Ds \vu) : \Ds \vuB } \dt + \int_0^\tau \intO{ \vr \vc{g} \cdot (\vu - \vu_b) }\dt. 
\end{split}
\]
However, a short inspection of the existence proof in \cite{ChoNovYan} reveals that \eqref{M3} can be established if suitable estimates 
yielding equi--integrability of the pressure potential $P$ are available at the level of approximate solutions. Seeing that 
\[
P(\vr) = P(\underline{\vr}) + (\vr - \underline{\vr} ) \left( P'(\underline{\vr}) - \frac{p(\underline{\vr})}{ 
\underline{\vr}} \right) + \vr \int_{\underline{\vr}}^{\vr} \frac{p(z)}{z^2} \ \D z \ 
\mbox{for any}\ 0 < \underline{\vr} < \vr < \Ov{\vr},
\]
we get
\[
P(\vr) \leq P(\underline{\vr}) + \Ov{\vr} \left| P'(\underline{\vr}) - \frac{p(\underline{\vr})}{ 
\underline{\vr}} \right| + \frac{ \Ov{\vr} }{ \underline{\vr}^2 } (\Ov{\vr} - \underline{\vr}) p(\vr)  
\]
yielding
\begin{equation} \label{mM3}
\limsup_{\vr \to \Ov{\vr} - } \frac{P(\vr)}{p(\vr)} = 0.
\end{equation}
Thus the pressure potential $P$ is dominated by the pressure $p$ near the singular value $\Ov{\vr}$. In particular, the bounds 
ensuring $L^1-$integrability of $p$ ensure equi--integrability of $P(\vr)$.

\begin{Remark}[Energy inequality] \label{RM1}

It is easy to see that the energy inequality \eqref{M3} is independent of the specific extension of the boundary velocity $\vu_b$ inside 
$\Omega$. Indeed, plugging $\psi (\vu^1_b - \vu^2_b)$, 
\[
\psi \in C^1_c(0, \infty),\ \vu^1_b|_{\partial \Omega} = \vu^2_b|_{\partial \Omega},
\]
as a test function in \eqref{M2} and $ \psi \frac12 (|\vu^1_b|^2 - |\vu^2_b|^2)$  as a test function in \eqref{M1}
we get 
\[
\begin{split}
\int_0^\infty &\intO{ \Big[ \vr \vu \cdot (\vu^1_b - \vu^2_b) \partial_t \psi +  \frac12 \vr  (|\vu^1_b|^2 - |\vu^2_b|^2) \partial_t \psi + \vr \vu \otimes \vu : \Grad (\vu^1_b - \vu^2_b) 
\psi
+ p(\vr) \Div (\vu^1_b - \vu^2_b) \psi \Big] } \dt\\ &= \int_0^\infty \intO{ \Big[ \mathbb{S}(\Ds \vu) : \Grad (\vu^1_b - \vu^2_b)  
- \vr \vc{g} \cdot (\vu^1_b - \vu^2_b)  -\vr \vu \cdot \Grad (|\vu^1_b|^2 - |\vu^2_b|^2)  \Big] \psi } \dt.
\end{split}
\]
\end{Remark}

\subsection{Geometry of the physical space}

We suppose that $\Omega \subset R^d$ be a bounded domain with $C^\infty$ boundary,  
\begin{equation} \label{EE1}
\partial \Omega = \cup_{i=0}^n \Gamma_i,  \ \Gamma_i \cap \Gamma_j = \emptyset,
\end{equation}
where $\Gamma_i$ are connected components of $\partial \Omega$, and $\Gamma_0$ is the boundary of the unbounded complement 
of $\Omega$ in $R^d$. We suppose
\begin{equation} \label{EE2}
\begin{split}
0 < \vr_b &< \Ov{\vr} \ \mbox{on}\ \Gamma_{\rm in},\\ 
\int_{{\Gamma_i}} \vu_b \cdot \vc{n} \ \D S_x &= 0,\ i=1, \dots, n,  \ 
\int_{\Gamma_0} \vu_b \cdot \vc{n}\ \D S_x > 0.
\end{split}
\end{equation}

\subsection{Main result}

Having collected the necessary preliminary material we are ready to state the main result. 

\begin{Theorem}[Existence of bounded absorbing set] \label{TM1}

Let $\Omega \subset R^d$, $d=2,3$ be a bounded domain of class $C^\infty$. Suppose that the boundary data $\vr_b$, $\vu_b$ satisfy \eqref{EE2} and that 
$\vc{g} \in L^\infty(\Omega; R^d)$.

Then there is $\mathcal{E}_\infty$ such that 
\begin{equation} \label{mM1}
\limsup_{t \to \infty} \intO{ \left[ \frac{1}{2} \vr |\vu|^2 + P(\vr) \right] (t, \cdot) } \leq \mathcal{E}_\infty,
\end{equation}
\begin{equation} \label{mM2}
\limsup_{T \to \infty} \int_T^{T+1} \left[ \left\| \vu \right\|_{W^{1,2}(\Omega; R^d)}^2 + \intO{ p(\vr) } \right] 
\dt \leq \mathcal{E}_\infty
\end{equation}
for any global--in--time solution $(\vr, \vu)$ of the Navier--Stokes system \eqref{i1}--\eqref{i3} specified in Definition \ref{MD1}.

\end{Theorem}

As observed in \eqref{mM3},
the pressure potential $P$ is dominated by the pressure $p$ near the singular value $\Ov{\vr}$ but not vice versa. 
$P$ may even stay bounded while $p$ blows up for $\vr \nearrow \Ov{\vr}$.

The rest of the paper is devoted to the proof of Theorem \ref{TM1}.

\section{Energy and pressure estimates}
\label{B}

In this section we collect the energy and pressure estimates available for the weak solutions of the Navier--Stokes system. 

\subsection{Energy inequality}

We start by rewriting the energy inequality \eqref{M3} in the form:

\begin{equation} \label{M4}
\begin{split}
&\left[ \intO{\left[ \frac{1}{2} \vr |\vu - \vuB|^2 + P(\vr) \right] } \right]_{t = T}^{t = T+\tau}  +
\int_T^{T + \tau} \intO{ \left( \mathbb{S}(\Ds \vu) - \mathbb{S}(\Ds \vuB) \right) : ( \Ds \vu - \Ds \vuB) } \dt \\
&+
\int_T^{T + \tau}  \int_{\Gamma_{\rm in}} P(\vr_b)  \vuB \cdot \vc{n} \ \D S_x \dt
\\	
&\leq -
\int_T^{T + \tau} \intO{ \left[ \vr (\vu - \vuB) \otimes (\vu - \vuB) + p(\vr) \mathbb{I} \right]  :  \Grad \vuB } \dt \\ & + \int_T^{T + \tau} \intO{ \vr \vc{g} \cdot (\vu - \vu_b) }\dt
- \int_T^{T + \tau} \intO{ \mathbb{S}(\Ds \vuB)  : ( \Ds \vu - \Ds \vuB) } \dt\\
&- \int_T^{T + \tau} \intO{  \vr \vuB \otimes (\vu - \vuB)  :  \Grad \vuB } \dt.
\end{split}
\end{equation}
Due to convexity of the energy functional 
\[
\frac{1}{2} \vr |\vu - \vuB|^2 + P(\vr) = 
\frac{|\vm|^2}{\vr} - \vm \cdot \vuB +  \vr |\vuB|^2 + P(\vr)
\]
with respect to the conservative variables $(\vr, \vm)$ that are weakly continuous in time, the energy inequality \eqref{M4} 
holds for a.a. $T > 0$ and all $\tau \geq 0$.

Consequently, under the hypotheses of Theorem \ref{TM1}, we may use Korn--Poincar\' e inequality and boundedness 
of $\vr$ to deduce 
\begin{equation} \label{S2}
\begin{split}
&\left[ \intO{\left[ \frac{1}{2} \vr |\vu - \vuB|^2 + P(\vr) \right] } \right]_{t = T}^{t = T+\tau}  +
\frac{\mu}{2} \int_T^{T + \tau} \| \vu - \vuB \|_{W^{1,2}_0(\Omega; R^d)}^2  \dt \\
&\leq -
\int_T^{T + \tau} \intO{ \left[ \vr (\vu - \vuB) \otimes (\vu - \vuB) + p(\vr) \mathbb{I} \right]  :  \Grad \vuB } \dt
+ \tau \omega\left(\Ov{\vr}, \| \vc{g} \|_{L^\infty}, \| \vuB \|_{W^{1, \infty}} \right),
\end{split}
\end{equation}
where the symbol $\omega$ denotes a generic function that is bounded for bounded arguments.
Note carefully that \eqref{S2} holds for \emph{any} extension of $\vuB$ inside $\Omega$. 

\subsection{Pressure estimates}

To derive estimates of the pressure $p$, we recall the so--called Bogovskii operator: 
\[
\begin{split}
\mathcal{B} : L^q_0 (\Omega) &\equiv \left\{ f \in L^q(\Omega) \ \Big| \ \intO{ f } = 0 \right\} 
\to W^{1,q}_0(\Omega, R^d),\ 1 < q < \infty,\\ 
\Div \mathcal{B}[f] &= f ,
\end{split}
\]  
see e.g. Galdi \cite[Chapter 3]{GALN}. 

For $\Phi \in L^q_0 (\Omega)$, we consider $\psi \mathcal{B}[\Phi]$, $\psi \in C^1_c(0, \infty)$ as a test function in the 
momentum equation \eqref{M2}: 
\begin{equation} \label{B1}
\begin{split}
\int_0^{\infty} \psi &\intO{ p(\vr) \Phi  } \dt =
- \int_0^{\infty} \partial_t \psi \intO{  \vr \vu \cdot \mathcal{B}[\Phi] } \dt\\
&- \int_0^{\infty} \psi \intO{  \vr \vu \otimes \vu : \Grad \mathcal{B}[\Phi] } \dt\\
&+ \int_0^{\infty} \psi \intO{  \mathbb{S} ( \Ds \vu ) : \Grad \mathcal{B} \left[ \Phi \right] } \dt - \int_0^{\infty} \psi \intO{  \vr \vc{g} \cdot \mathcal{B} \left[ \Phi \right] } \dt
\end{split}
\end{equation}
As the momentum $\vr \vu$ is weakly continuous in time, we may infer that 
\begin{equation} \label{B2}
\begin{split}
\int_T^{T + \tau} &\intO{ p(\vr) \Phi } \dt =
\left[ \intO{ \vr \vu \cdot \mathcal{B} \left[ \Phi \right] }  \right]_{t=T}^{t = T + \tau} \\
&- \int_T^{T + \tau} \intO{  \vr \vu \otimes \vu : \Grad \mathcal{B} \left[ \Phi \right] } \dt
+ \int_T^{T + \tau} \intO{  \mathbb{S} ( \Ds \vu ) : \Grad \mathcal{B} \left[\Phi \right] } \dt \\
&- \int_T^{T + \tau} \intO{  \vr \vc{g} \cdot \mathcal{B} \left[ \Phi \right] } \dt
\end{split}
\end{equation}
for any $T > 0$, $\tau \geq 0$.

\section{Extending boundary vector fields}
\label{E} 

We construct a suitable extension of the boundary vector field $\vuB$ inside $\Omega$.

\subsection{First extension}

First observe that the outer component of the boundary
$\Gamma_0$ contains at least one \emph{extremal point} $x_0 \in \Gamma_0$, specifically,  
\[
\Ov{\Omega} \cap \tau_{x_0} = x_0,\ \mbox{where}\ \tau_{x_0} 
\ \mbox{denotes the tangent plane to}\ \partial \Omega \ \mbox{at}\ x_0.
\]
Without loss of generality, we may assume that
\[ \Omega \subseteq \{x: x^1 <x_0^1\} \mbox{ and }
x_0 = [x_0^1,0,\dots,0],\ \tau_{x_0} = x_0 + R^{d-1}.
\]
Consider a function 
\[
\chi(z) = \left\{ \begin{array}{l} \ 0 \ \mbox{if}\ z \leq 0,\\ 
\chi'(z) > 0 \ \mbox{for}\ z > 0 \end{array} \right., 
\]
together with a vector field
\[
\vc{v}^0_b(x) = \lambda \left[ \chi (x^1 - x_0^1 + \delta), 0 , \dots, 0 \right].  
\]
It is easy to check that 
\[
\Ds \vc{v}^0_b = \begin{bmatrix} \lambda \chi' (x^1 - x_0^1 + \delta)  & 0 & 0 \\ 
0 & 0 &0 \\
0 & 0 & 0
\end{bmatrix},\ \Div \vc{v}^0_b = \lambda \chi' (x^1 - x_0^1 + \delta).
\]
Next, we choose $\delta > 0$ small enough so that 
\[
\vc{v}^0_b |_{\Gamma_i } = 0, \ i =1, \dots, n,\ \vc{v}^0_b|_{\Gamma _0} \ne 0, 
\]
and then $\lambda > 0$ large enough so that 
\[
\int_{\Gamma_0} \vc{v}^0_b \cdot \vc{n} \ \D S_x = \int_{\Gamma_0} \vu_b \cdot \vc{n}\ \D S_x > 0.
\]
Finally, we set 
\[
\vuB = \vc{w}_b + \vc{v}^0_b, 
\] 
where $\vc{w}_b = \vu_b - \vc{v}^0_b$, and
\begin{equation} \label{EE3}
\begin{split}
\int_{\Gamma_i} \vc{w}_b \cdot \vc{n} \ \D S_x &= 0 \ \mbox{for all}\ i = 0,1,\dots, n,\\
\Ds \vc{v}^0_b &\geq 0,\  
\ \mbox{there is an open set}\ B \subset \Omega, \ |B| > 0,\ 
\inf_B (\Div \vc{v}^0_b) \geq D > 0.
\end{split}
\end{equation} 

\subsection{Extension of the vector field $\vc{w}_b$}

As $\vc{w}_b$ satisfies \eqref{EE3}, we can use  lemma by Galdi \cite[Lemma IX.4.1]{GALN}
(see also the decomposition theorem of Kozono and Yanagisawa \cite[Proposition 1]{KozYan})
to write 
\[
\vc{w}_b = {\bf curl} \ \vc{z}_b \ \mbox{in}\ \Ov{\Omega}, 
\]
where $\vc{z}_b$ is smooth provided $\vc{w}_b$ is smooth. 

Next, we report the following result by Galdi \cite[Lemma III.6.1,Lemma III.6.2]{GALN}: For each $\ep > 0$, there exists a function 
$d_\ep \in C^\infty(\Ov{\Omega})$ enjoying the following properties:
\begin{itemize}
\item 
\begin{equation} \label{EE4}
|d_\ep | \leq 1,\ d_\ep (x) \equiv 1 \ \mbox{for all}\ x \ \mbox{in an open neighborhood of}\ \partial \Omega;
\end{equation}
\item 
\begin{equation} \label{EE5}
d_\ep (x) \equiv 0 \ \mbox{whenever}\ {\rm dist}[x, \partial \Omega] > \ep;
\end{equation}
\item
\begin{equation} \label{EE6}
|D^\alpha_x d_\ep (x) | \leq c \frac{\ep}{{\rm dist}^{|\alpha|} [x, \partial \Omega]},\ |\alpha| = 1,2 ,\ 0 < \ep < 1,\ 
x \in \Omega;
\end{equation}
where the constant is independent of $\ep$.

\end{itemize}

Finally, we extend $\vuB$ as 
\begin{equation} \label{EE7}
\vu_b^\ep = {\bf curl} (d_\ep \vc{z}_b ) + \vc{v}^0_b.
\end{equation}

\section{Bounded absorbing set -- proof of the main result}
\label{S}

Going back to the inequality \eqref{S2} and using the decomposition \eqref{EE7}, the integral on the right--hand side 
may be handled as
 \[
\begin{split} 
-
\int_T^{T + \tau} &\intO{ \left[ \vr (\vu - \vuB^\ep) \otimes (\vu - \vuB^\ep) + p(\vr) \mathbb{I} \right]  :  \Grad \vuB^\ep } \dt
\\
&=-  \int_T^{T + \tau} \intO{ \left[ \vr (\vu - \vuB^\ep) \otimes (\vu - \vuB^\ep) + p(\vr) \mathbb{I} \right]  :  \Grad \vc{v}^0_b } \dt
\\
&- \int_T^{T + \tau} \intO{ \left[ \vr (\vu - \vuB^\ep) \otimes (\vu - \vuB^\ep) + p(\vr) \mathbb{I} \right]  :  \Grad {\bf curl}
(d_\ep \vc{z}_b ) } \dt\\
 &\leq - D \int_T^{T + \tau} \int_B { p(\vr) } \dx \dt + \Ov{\vr} \int_T^{T + \tau} \int_{\Omega} \frac{|\vu - \vuB^\ep |^2}{
{\rm dist}^2(x, \partial \Omega)} {\rm dist}^2(x, \partial \Omega) 
|\Grad {\bf curl} (d_\ep \vc{z}_b) | \dx \dt,  
\end{split}
\]
where the open set $B$ has been identified in \eqref{EE3}.

Recalling Hardy--Sobolev inequality 
\[
\intO{ \frac{|\vu - \vuB^\ep |^2}{{\rm dist^2(x, \partial \Omega)}}}   
\aleq \| \vu - \vuB^\ep \|_{W^{1,2}_0(\Omega; R^d)}^2,
\]
and the properties of $d_\ep$, specifically \eqref{EE6}, we may 
fix $\ep > 0$ small enough so that  inequality \eqref{S2} reduces to 
\begin{equation} \label{S3}
\begin{split}
&\left[ \intO{\left[ \frac{1}{2} \vr |\vu - \vuB^\ep |^2 + P(\vr) \right] } \right]_{t = T}^{t = T+\tau}  +
\frac{\mu}{4} \int_T^{T + \tau} \| \vu - \vuB^\ep \|_{W^{1,2}_0(\Omega; R^d)}^2  \dt + 
D \int_T^{T + \tau} \int_B { p(\vr) } \dx \dt
\\
&\leq 
\tau \omega\left(\Ov{\vr}, \| \vc{g} \|_{L^\infty}, \| \vuB^\ep \|_{W^{1, \infty}} \right)
\end{split}
\end{equation}
for any $\ep > 0$.  

Finally, evoking the pressure estimates \eqref{B2}, we get, following the same line of arguments,
\begin{equation} \label{S7}
\begin{split}
\int_T^{T + \tau} &\intO{ p(\vr) \Phi } \dt \leq
\left[ \intO{ \vr (\vu - \vuB^\ep) \cdot \mathcal{B} \left[ \Phi \right] }  \right]_{t=T}^{t = T + \tau} \\ &+ c(\Ov{\vr}, \Phi) \int_T^{T + \tau} \| \vu - \vuB^\ep \|^2_{W^{1,2}_0(\Omega; R^d)} \dt + \tau \omega\left(\Ov{\vr}, \| \vuB^\ep \|_{W^{1, \infty}}, \Phi \right)
\end{split}
\end{equation} 
for any 
\[
\Phi \in L^q(\Omega), \ \intO{ \Phi } = 0.
\]

Now, in \eqref{S7} we choose $\Phi\in L^\infty_0(\Omega)$ such that $\Phi|_{\Omega \setminus B} = 1$ to obtain 
\begin{equation} \label{S9}
\begin{split}
\int_T^{T + \tau} &\int_{\Omega \setminus B} p(\vr)  \ \dx \dt \leq
\left[ \intO{ \vr (\vu - \vuB^\ep) \cdot \mathcal{B} \left[ \Phi \right] }  \right]_{t=T}^{t = T + \tau} \\ &+ c(\Ov{\vr}, \Phi) \int_T^{T + \tau} \| \vu - \vuB^\ep \|^2_{W^{1,2}_0(\Omega; R^d)} \dt + \tau \omega\left(\Ov{\vr}, \| \vuB^\ep \|_{W^{1, \infty}}, \Phi \right)
+ \int_T^{T + \tau} \int_{B} p(\vr) |\Phi| \ \dx \dt.
\end{split}
\end{equation} 
Thus, going back to \eqref{S3}, for all $\delta > 0$ small enough,  
\begin{equation} \label{S12}
\begin{split}
&\left[ \intO{\left[ \frac{1}{2} \vr |\vu - \vuB^\ep |^2 + P(\vr) \right] } \right]_{t = T}^{t = T+\tau}  +
\frac{\mu}{8} \int_T^{T + \tau} \| \vu - \vuB^\ep \|_{W^{1,2}_0(\Omega; R^d)}^2  \dt + \delta
\int_T^{T + \tau} \intO{ p(\vr)} \dt \\
&\leq 2 \delta \left[ \intO{ \vr (\vu - \vuB^\ep ) \cdot \mathcal{B} \left[ \Phi \right] }  \right]_{t=T}^{t = T + \tau} + \tau \omega\left(\Ov{\vr}, \| \vc{g} \|_{L^\infty}, \| \vuB^\ep \|_{W^{1, \infty}}, \Phi \right)
\end{split}
\end{equation}

Seeing that $P(\vr) \aleq p(\vr)$ for $\vr \nearrow\Ov{\vr}$ (see \eqref{mM3}) we obtain the desired conclusion.

\section{Concluding remarks}
\label{C}

As we have seen, validity of Theorem \ref{TM1} depends essentially on the hypothesis \eqref{EE2} strongly reminiscent of the 
so--called Leray's condition, see e.g. Galdi \cite{Gald1}. A short inspection of the proof reveals that Theorem 
\eqref{TM1} remains valid if 
\begin{equation} \label{C1}
\vu_b = {\bf curl} \ \vc{z}_b + \vu^0_b \ \mbox{in}\ \Ov{\Omega},  
\end{equation}
where 
\begin{equation} \label{C2}
\Ds \vu^0_b \geq 0, \ \inf_B (\Div \vu^0_b) > 0 \ \mbox{for some non--empty open set}\ B \subset \Omega.
\end{equation}
Of course, Theorem \ref{TM1} holds under the general hypothesis \eqref{i9} as soon as $\Omega$ is simply connected.

Note that the conclusion of Theorem \ref{TM1} may fail in the case 
\[
\int_{\partial \Omega} \vu_b \cdot \vc{n} \ \D S_x < 0.
\]
Indeed, for $\Gamma_{{\rm in}} = \partial \Omega$, $\vr_b > 0$, we get 
\[
\intO{ \vr(t, \cdot) } \to \infty \ \mbox{as}\ t \to \infty 
\]
in contrast with $\vr \leq \Ov{\vr}$.

In the case of the isentropic EOS $p(\vr) = a \vr^\gamma$, the arguments are hampered by the lack of control over the convective term 
\[
- \vr (\vu - \vuB) \otimes (\vu - \vuB) : \Grad \vu_b
\]
on the right--hand side of the energy inequality \eqref{S2}. The easy case $\Ds \vu_b \geq 0$ is treated in \cite{FanFeiHof}.

\subsection{Total mass conservation, tangential boundary velocity}

An interesting situation occurs when the boundary field $\vuB$ is tangential to $\partial \Omega$. Obviously,
\[
\int_{\partial \Omega} \vu_b \cdot \vc{n} \ \D S_x = 0, 
\]
in this case and Theorem \ref{TM1} does not apply. Note that such a scenario causes Taylor instability in turbulence theory, where the fluid is excited by one or more rotating bodies, see e.g. Davidson \cite{DAVI}.

First observe that, in accordance with Kozono and Yanagisawa \cite{KozYan}, the vector field $\vuB$ admits an extension 
\[
\vu_b^\ep = {\bf curl}(d_\ep \vc{z}_b) 
\]
as in \eqref{EE7}. Consequently, following the line of arguments used in Section \ref{S} we arrive at the inequality 
\begin{equation} \label{C3}
\begin{split}
\left[ \intO{\left[ \frac{1}{2} \vr |\vu - \vuB^\ep |^2 + P(\vr) \right] } \right]_{t = T}^{t = T+\tau}  +
\frac{\mu}{4} \int_T^{T + \tau} \| \vu - \vuB^\ep \|_{W^{1,2}_0(\Omega; R^d)}^2  \dt \leq 
\tau \omega\left(\Ov{\vr}, \| \vc{g} \|_{L^\infty}, \| \vuB^\ep \|_{W^{1, \infty}} \right).
\end{split}
\end{equation}

Next, we derive refined pressure estimates based on the fact that $\vu \cdot \vc{n}|_{\partial \Omega} = 0$. We consider  
\[
\bfphi = \psi \mathcal{B} \left[ \vr - \frac{1}{|\Omega|} \intO{ \vr } \right],\ \psi \in C^1_c(0, \infty)  
\]
as a test function in the momentum balance \eqref{M2}. After a straightforward manipulation, we deduce 
\begin{equation} \label{C4}
\begin{split}
\int_T^{T + \tau} &\intO{ p(\vr) \left[ \vr - \frac{1}{|\Omega
|} \intO{ \vr } \right]  } \dt =
\left[ \intO{ \vr \vu \cdot \mathcal{B} \left[ \vr - \frac{1}{|\Omega
|} \intO{ \vr } \right] }  \right]_{t=T}^{t = T + \tau} \\
&+ \int_T^{T + \tau} \intO{  \vr \vu \cdot \mathcal{B} [ \Div (\vr \vu) ] } \dt\\
&- \int_T^{T + \tau} \intO{  \vr \vu \otimes \vu : \Grad \mathcal{B} \left[ \vr - \frac{1}{|\Omega|} \intO{ \vr } \right] } \dt\\
&+ \int_T^{T + \tau} \intO{  \mathbb{S} ( \Ds \vu ) : \Grad \mathcal{B} \left[ \vr - \frac{1}{|\Omega|} \intO{ \vr } \right] } \dt \\
&- \int_T^{T + \tau} \intO{  \vr \vc{g} \cdot \mathcal{B} \left[ \vr - \frac{1}{|\Omega|} \intO{ \vr } \right] } \dt
\end{split}
\end{equation}
At this stage, we recall that $\mathcal{B}: L^q_0 \to W^{1,q}_0$ is a bounded linear operator for any $1 < q < \infty$, and, in addition, 
\[
\left\| \mathcal{B} [\Div \vc{v} ] \right\|_{L^r(\Omega; R^d)} \aleq \| \vc{v} \|_{L^r(\Omega; R^d)} 
\ 1 < r < \infty \ \mbox{if}\ \Div \vc{v} \in L^q(\Omega), \ \vc{v} \cdot \vc{n}|_{\partial \Omega} = 0, 
\]
see Geissert, Heck, and Hieber \cite{GEHEHI}. Consequently, all the integrals on the right--hand side may be controlled by the dissipation term in \eqref{C3} exactly as in Section \ref{S}. 

Finally, as the total mass $M$ is a constant of motion, and, consequently 
\[
\frac{M}{|\Omega|} = \frac{1}{|\Omega|} \intO{ \vr } < \Ov{\vr},
\]
we get
\[  
\intO{ p(\vr) \left[ \vr - \frac{1}{|\Omega
|} \intO{ \vr } \right]  } \ageq \intO{ p(\vr) } - c(M),
\]
which yields the following result.

\begin{Theorem}[Impermeable boundary] \label{TC1}

Let $\Omega \subset R^d$, $d=2,3$ be a bounded domain of class $C^\infty$. Suppose that  
$\vc{g} \in L^\infty(\Omega; R^d)$, and that the boundary velocity satisfies 
\[
\vuB \cdot \vc{n}|_{\partial \Omega} = 0.
\]

Then there is $\mathcal{E}_\infty$ such that 
\[
\limsup_{t \to \infty} \intO{ \left[ \frac{1}{2} \vr |\vu|^2 + P(\vr) \right] (t, \cdot) } \leq \mathcal{E}_\infty,
\]
\[
\limsup_{T \to \infty} \int_T^{T+1} \left[ \left\| \vu \right\|_{W^{1,2}(\Omega; R^d)}^2 + \intO{ p(\vr) } \right] 
\dt \leq \mathcal{E}_\infty
\]
for any global--in--time solution $(\vr, \vu)$ of the Navier--Stokes system \eqref{i1}--\eqref{i3} specified in Definition \ref{MD1}.

\end{Theorem}
\def\cprime{$'$} \def\ocirc#1{\ifmmode\setbox0=\hbox{$#1$}\dimen0=\ht0
  \advance\dimen0 by1pt\rlap{\hbox to\wd0{\hss\raise\dimen0
  \hbox{\hskip.2em$\scriptscriptstyle\circ$}\hss}}#1\else {\accent"17 #1}\fi}


\end{document}